\documentclass[12pt]{article}

\def\ra{\rightarrow}

\def\ss{\subseteq}

\def\d{\delta}

 \def\HollowBox #1#2{{\dimen0=#1 \advance\dimen0 by -#2       
       \dimen1=#1 \advance\dimen1 by #2                       
        \vrule height #1 depth #2 width #2                    
        \vrule height 0pt depth #2 width #1                   
        \llap{\vrule height #1 depth -\dimen0 width \dimen1}%
       \hskip -#2                                             
       \vrule height #1 depth #2 width #2}}                   
 \def\BoxOpTwo{\mathord{\HollowBox{6pt}{.4pt}}\;}             

\def\endpf{\hfill $\BoxOpTwo$}

\font\teneufm=eufm10
\font\seveneufm=eufm7
\font\fiveeufm=eufm5
\newfam\eufmfam
\textfont\eufmfam=\teneufm
\scriptfont\eufmfam=\seveneufm
\scriptscriptfont\eufmfam=\fiveeufm

\newfam\msbfam
\font\tenmsb=msbm10 scaled \magstep1 \textfont\msbfam=\tenmsb
\font\sevenmsb=msbm7 scaled \magstep1 \scriptfont\msbfam=\sevenmsb
\font\fivemsb=msbm5  scaled \magstep1  \scriptscriptfont\msbfam=\fivemsb
\def\Bbb{\fam\msbfam \tenmsb}

\def\RR{{\Bbb R}}

\newtheorem{theorem}{Theorem}

\newtheorem{proposition}[theorem]{Proposition}
\newtheorem{lemma}[theorem]{Lemma}

\newtheorem{definition}{Definition}

\makeindex

\usepackage{graphicx}

\usepackage{amsmath}

\begin{document}

\begin{center}
\huge \bf
Some New Thoughts on Maximal Functions and Poisson Integrals
\end{center}
\vspace*{.12in}

\begin{center}
\large Steven G. Krantz\footnote{Author supported in part
by the National Science Foundation and by the Dean of the Graduate
School at Washington University.}\footnote{{\bf MR Classification
Numbers:} 42B25, 42B99, 31B05, 31B10.}
\end{center}
\vspace*{.15in}

\begin{center}
\today
\end{center}
\vspace*{.2in}

\begin{quotation}
{\bf Abstract:} \sl
We study Wiener-type covering lemmas, Hardy-Littlewood-type maximal
functions, and convergence theorems on metric spacs.   Later we specialize
down to a result for the Poisson integral.  We show
that, in a suitably general setting, these three phenomena are essentially
logically equivalent.  Along the way we discuss some useful estimates
for the Poisson kernel.
\end{quotation}
\vspace*{.25in}

\setcounter{section}{-1}

\section{Introduction}

A standard paradigm in harmonic analysis---on Euclidean space,
or on a space of homogeneous type---is that if one can prove a
covering lemma of Wiener type then one can prove a weak-type
estimate for a suitable maximal function. And this, in turn,
will imply pointwise convergence results for certain
convolution operators. The celebrated ``limits of sequences of
operators'' theorem of Stein [STE1] fleshes out this picture
by showing that, under suitable hypotheses, a pointwise
convergence result implies a maximal function estimate. One of
our goals in this paper is to complete this logical scenario.
In particular, we wish to show that if one has pointwise
convergence of integral operators, then one can derive a
covering lemma. Thus, in effect, the three key ideas being
discussed here are logically equivalent.

We begin our investigations with a few words about the size of the
Poisson kernel.

\section{Estimates on the Poisson Kernel}

Of course the Poisson kernel on classical domains is well known.
For example,
\begin{itemize}
\item The Poisson kernel of the disc $D \ss \RR^2$ is
$$
P_D(z,\zeta) = \frac{1}{2\pi} \cdot \frac{1 - |z|^2}{|z - \zeta|^2} \, ,
$$
where $z \in D$ and $\zeta \in \partial D$.
\item The Poisson kernel for the upper halfplane
$$
U^2 = \{(x_1, x_2) \in \RR^2: x_2 > 0\}
$$
is given by
$$
P_{U^2}(x,t) = \frac{1}{\pi} \cdot \frac{x_2}{(x_1 - t)^2 + x_2^2} \, ,
$$
where $x = (x_1, x_2) \in U^2$ and $t \in \RR = \partial U^2$.
\item The Poisson kernel for the unit ball $B \ss \RR^N$
is given by
$$
P_B(x,t) = \frac{\Gamma(N/2)}{2 \pi^{N/2}} \cdot \frac{1 - |x|^2}{|x - t|^N} \, ,
$$
where $x = (x_1, x_2, \dots, x_N) \in B$ and $t = (t_1, \dots, t_N) \in \partial B$.
Here $\Gamma$ is the classical gamma function.
\item The Poisson kernel for the upper halfspace $U^{N+1} \equiv
\{x = (x_1,\dots, x_{N+1}) \in \RR^{N+1}: x_{N+1} > 0\}$
(with $x = (x_1,\dots, x_{N+1}) \equiv (x', x_{N+1})$) is given by
$$
P_{U^{N+1}}(x,t) = c_N \frac{x_{N+1}}{([x' - t]^2 + x_{N+1}^2)^{[N+1]/2}} \,
$$
where $x = (x_1, x_2, \dots, x_{N+1}) \in U^{N+1}$, $t = (t_1,t_2,\dots, t_N) \in \RR^N = \partial U^{N+1}$, and
$$
c_N = \frac{\Gamma([N+1]/2)}{\pi^{[N+1]/2}} \, .
$$
\end{itemize}

We say that $\Omega \ss \RR^N$ is a {\it domain} if it is a connected,
open set.  Of course there is no hope of producing an actual formula for
the Poisson kernel of an arbitrary domain in $\RR^N$.  But there is definite interest in obtaining 
an {\it asymptotic} formula for the Poisson kernel on a general domain.  The standard
asymptotic (see [STE3] or [KRA1]) is
$$
P_\Omega(x,t) \approx \frac{\delta(x)}{|x - t|^N} 	\eqno (*)
$$
for $x \in \Omega$ and $t \in \partial \Omega$.
Here $\delta(x) \equiv \delta_\Omega(x)$ is the distance
from $x \in \Omega$ to $\partial \Omega$.  This estimate, together
with analogous estimates for the derivatives of $P_\Omega$,
suffices for most applications.  Stein states this result in [STE3], but
does not prove it.  The reference that he gives is also incomplete in this
matter.  Perhaps the first complete proof to appear in print can be found in [KRA1].
A more recent, and more efficient, proof appears in [KRA2].  In fact a very
natural way to approach the matter is to fix a point $x$ in $\Omega$ near the
boundary and to consider a small, topologically trivial subdomain $\Omega' \ss \Omega$ 
which contains $x$ and which shares a piece $W$ of boundary with $\Omega$.  It follows
from basic estimates on the Green's function (see, for instance, [APF]) that
the Poisson kernel of $\Omega$ for this fixed point $x$ and for $t \in W$ is comparable
to the Poisson kernel of $\Omega'$ for that same $x$ and $t$ (the constants of comparison, of
course, depend on $\Omega$ and $\Omega'$).   Now of $\Omega'$ may be mapped 
diffeomorphically to the unit ball by a mapping
$\Phi$.  And the Poisson kernel $P$ of the unit ball may be pulled back to
$\Omega'$ by means of pseudodifferential operator theory.  One obtains immediately
that
$$
P_{\Omega'}(x,y) \approx \frac{\delta_{\Omega'}(x)}{|x - y|^N} + {\cal E}(x,y) \, ,
$$
where ${\cal E}(x,y)$ is an error term of lower order.  The result follows.

We will make good use of the approximation $(*)$ in the sequel.

\section{Principal Results}

We begin by enunciating the results that are known.  We do so in the setting of an
arbitary metric space ({\it not} a space of homogeneous type).  

We first need some definitions.	 In what follows, let $(X, \rho)$ be a metric
space.   The balls in this metric space will be denoted by
$$
B(x,r) \equiv \{t \in X: \rho(x,t) < r\} \, .
$$
We take it that the metric space is equipped with a Borel regular measure (i.e., a Radon measure) $\mu$, and that
the measure of each ball is positive and finite.  But we do {\it not} assume
that $(X, \rho, \mu)$ is a space of homogeneous type in the sense of [COW].
Sometimes, for convenience, we will denote the $\mu$-measure
of a ball $B(x,r)$ by $|B(x,r)|$ or, more generally, the $\mu$-measure
of a measurable set $S$ by $|S|$.

\begin{definition} \rm
We say that the balls $B(x,r)$
satisfy a {\it Wiener-type covering lemma} if, whenever $K \ss X$ is compact
and $\{B_\alpha\}_{\alpha \in A}$ is a covering of $K$ by balls, then
there is a subcollection $\{B_{\alpha_j}\}_{j=1}^K$ such that
\begin{enumerate}
\item[{\bf (a)}]  The $\{B_{\alpha_j}\}$ are pairwise disjoint.
\item[{\bf (b)}]  The dilated balls $\{3B_{\alpha_j}\}$, where
$3B(x,r) \equiv B(x, 3r)$, cover $K$.
\end{enumerate}
\end{definition}

\begin{definition} \rm
A real-valued function $f$ on an open set $U \ss X$ is said
to be of {\it weak type $\alpha$}, $\alpha > 0$, if there is
a $C > 0$ such that, for any $\lambda > 0$,
$$
\mu\{x \in U: |f(x)| > \lambda\} \leq \frac{C}{\lambda^\alpha} \, .
$$
\end{definition}
It is easy to see that any $L^\alpha$ function is in fact
weak type $\alpha$.

\begin{definition} \rm
A linear operator defined on $L^p(U)$, for some open set $U \ss X$, and
taking values in the real functions on some set $V \ss X$, is said
to be of {\it weak type $(p,p)$} if there is a constant $C > 0$ so that,
for any $\lambda > 0$, 
$$
\mu\{x \in V: |Tf(x)| > \lambda\} \leq C \cdot \left [\frac{\|f\|_{L^p(U)}}{\lambda} \right ]^p \, .  \eqno (\star)
$$
We say that the linear operator is of {\it restricted weak type} $(p,p)$ if
it satisfies the ``weak type $(p,p)$'' condition $(\star)$ 
when restricted to act on $f$ the characteristic function
of a set.  
\end{definition}

\begin{definition}   \rm
For a real-valued, locally integrable function $f$ defined on an open $U \ss X$, and $x \in U$, define
$$
M f(x) = \sup_{r > 0} \frac{1}{|B(x,r)|} \, \int_{B(x,r)} |f(t)| \, dt \, .
$$
We call $M$ the {\it Hardy-Littlewood maximal operator} or {\it Hardy-Littlewood maximal
function}.
\end{definition}

\begin{definition} \rm
Let ${\cal B} = \{B_\alpha\}_{\alpha \in A}$ be a covering of some set $E$ by balls.
We call ${\cal C} = \{C_\beta\}_{\beta \in B}$ a {\it refinement} of the covering ${\cal B}$
if each element $C_\beta$ of ${\cal C}$ is a subset of some $B_\alpha \in {\cal B}$ and
if $\cup_\beta C_\beta \supset E$.
\end{definition}

\begin{definition} \rm
We say that the {\it Lebesgue differentiation theorem} holds in $L^p(X)$ if,
for any $f \in L^p(X)$, and for almost every $x \in X$, it holds that
$$
\lim_{r \ra 0^+} \frac{1}{|B(x,r)|} \, \int_{B(x,r)} f(t) \, dt
$$
exists and equals $f$.
\end{definition}

Now we have

\begin{proposition} \sl
Assume that there is a constant $C > 0$ so that
$$
\mu(B(x,3r)) \leq C \cdot \mu(B(x,r) 
$$
for any $x \in X$ and any $r >0$.  Then the fact that metric balls $B(x,r)$ satisfy a Wiener-type covering lemma
implies that the Hardy-Littlewood maximal operator is weak type $(1,1)$.

The fact that the Hardy-Littlewood maximal operator is weak type $(1,1)$ implies
that the Lebesgue differentiation theorem holds in $L^1$.

Under suitable additional hypotheses, the fact that the
Lebesgue differentiation theorem holds for functions in $L^1$
implies that the Hardy-Littlewood maximal operator is weak
type $(1,1)$. 
\end{proposition}

These statements are all classical.  The third one is---in the classical setting
of the circle group, for instance---the content
of Stein's limits of sequences of operators theorem [STE1].  In our more
general setting of metric spaces we would want to apply instead Sawyer's generalization
of this result [SAW] or Nikisin's generalization of this result [GIL]; this
would require the averaging operators to satisfy some ergodic or mixing condition.

Our intention now is to round out the logic in this last
proposition. We shall show that the fact that the Lebesgue
differentiation theorem holding for functions in $L^1$ implies
that metric balls $B(x,r)$ satisfy a covering 
lemma.

We have intentionally omitted the phrase ``Wiener-type'' from the enunciation
of this last result because in fact we will need to formulate a new
covering lemma in order to make the logic work.  This new covering
lemma is of intrinsic interest, and we shall spend some time putting
it into context.  

\section{The New Covering Lemma}

In order for the logic to mesh properly, we need to formulate
a new type of covering lemma.  Then we need to show that
it does the same job as the old (Wiener-type) covering lemma.
The new lemma is this (we work still on the metric space $(X,\rho)$, which
is {\it not} necessarily a space of homogeneous type):

\begin{lemma} \sl
The Lebesgue differentiation theorem implies the following covering lemma:
Let $\{B_\alpha\}_{\alpha \in A}$ be a covering of a compact set $K$ by
open metric balls.  Then there is a refinement $\{C_j\}$ of the covering so that
\begin{enumerate}
\item[{\bf (a)}]  The balls $C_j$ are pairwise disjoint.
\item[{\bf (b)}]  Each ball $C_j$ has the property that 
$$
\frac{|C_j \cap K|}{|C_j|} > \frac{1}{2} \, ;
$$
\item[{\bf (c)}]  We have the inequality 
$$
\sum_j |C_j| > \frac{1}{2} |K| \, ;
$$
\end{enumerate}
\end{lemma}
{\bf Proof:}  Fix a compact set $K \ss X$.  For almost every $x \in K$,
the Lebesgue differentiation theorem tells us that there is an $r_x > 0$
so that, if $0 < s \leq r_x$ then
$$
|B(x,s) \cap K| > \frac{1}{2} |B(x,s)| \, .
$$
By Borel regularity, let $U$ be an open set that contains $K$ and so that $|U \setminus K| < (1/10) \cdot |K|$.
We may assume that each $r_x$ is so small that $B(x, r_x) \ss U$ for each $x \in K$.

Now if $x \in U \setminus K$ then take $r_x = \min(\hbox{dist}(x, K), \hbox{dist}(x, {}^c U)$.
We may apply the (proof of) the Vitali covering theorem now to extract a subcollection
$\{B(x_j, r_{x_j})$ so that the $B(x_j, r_{x_j})$ are pairwise disjoint
and so that $\sum_j |B(x_j, r_{x_j})| > (1/2) K$.  
\endpf 
\smallskip \\

\begin{proposition} \sl
The Lebesgue differentiation theorem implies that the 
Hardy-Littlewood maximal
operator is restricted weak type $(1,1)$.
\end{proposition}
\noindent {\bf Proof:}   
Fix a measurable set $E$ with $|E| < \infty$.
Now let $0 < \lambda < 1$ and set
$$
S_\lambda = \{x \in X: M \chi_E > \lambda\} \, .
$$
Let $K \ss S_\lambda$ be compact and satisfy $|K| > |S_\lambda|/2$.
For each $x \in K$ there is a ball $B_x$ such that
$$
\frac{1}{|B_x|} \, \int_{B_x} \chi_E (t) \, dt > \lambda \, .
$$

Now, as in the proof of Lemma 2, choose a refinement $\{B_{x_j}\}$ so that, for each $j$, 
$$
\frac{|B_{x_j} \cap K|}{B_{x_j}} > \frac{1}{2} 
$$
and so that
the $\{B_{x_j}\}$ are pairwise disjoint.  Finally we may assume
that $\sum_j |B_{x_j}| > |K|/2$, again by arguing as in the
Vitali covering lemma.  

Then
\begin{eqnarray*}
    |S_\lambda| & \leq & 2|K| \\
             & \leq & 4 \sum_j |B_{x_j}| \\
             & \leq & \frac{4}{\lambda} \cdot \sum_j |B_{x_j} \cap E | \\
	     & \leq & \frac{4|E|}{\lambda} \, .
\end{eqnarray*}
That is the estimate that we want.  So $M$ is restricted weak
type $(1,1)$.
\endpf 
\smallskip \\

Now the sum of what we have established thus far is this:
\smallskip \\
$$
\hbox{\bf (covering lemma)} \Longrightarrow \hbox{\bf (maximal function estimate)} 
$$
\vspace*{-.25in}
$$
\ \ \ \ \ \ \ \ \ \ \ \ \ \ \Longrightarrow
    \hbox{\bf (differentiation theorem)} \Longrightarrow \hbox{\bf (covering lemma)}
$$
$$
 \Longrightarrow
    \hbox{\bf (maximal function estimate)}
$$
\vspace*{.15in}

\noindent So we have closed the logical circle described at the beginning of this paper.

\section{More General Results}

In fact the Lebesgue differentiation theorem is merely a paradigm for the
type of result that is true in general.  We first introduce some notation
and terminology.

If $\Omega \ss \RR^N$ is a smoothly bounded domain, then we let $d\sigma$ denote
$(N-1)$-dimensional Hausdorff measure on $\partial \Omega$.  For $x \in \Omega$,
let $\delta(x)$ denote the distance of $x$ to $\partial \Omega$.  If $y \in \partial \Omega$
and $r > 0$ then let $\beta(y,r) = \{t \in \partial \Omega: |t - y| < r\}$.
If $y \in \partial \Omega$ and $\alpha > 0$, then let
$$
\Gamma_\alpha(y) = \{x \in \Omega: |x - y| < \alpha\cdot \delta(x)\} \, .
$$
								 
\begin{theorem} \sl
Let $\Omega$ be a smoothly bounded domain in $\RR^N$.  Let $P(x,t)$
be the Poisson kernel for $\Omega$, with $x \in \Omega$ and $t \in \partial \Omega$.
If $f \in L^1(\partial \Omega, d\sigma)$, then set
$$
u(x) = \int_{\partial \Omega} P(x,t) f(t) \, d\sigma(t) \, .
$$
Fix $\alpha > 0$.  The $\sigma$-almost everywhere existence of the boundary
limits
$$
\lim_{\Gamma_\alpha(y) \ni x \ra y} u(x) = f(y) 
$$
is logically equivalent to the restricted weak-type $(1,1)$ property of the maximal operator
$$
{\cal M} f(y) \equiv \sup_{r > 0} \frac{1}{\sigma(\beta(y,r)} \, \int_{\beta(y,r)} |f(t)| \, d\sigma(t) \, .
$$
\end{theorem}
{\bf Proof:}  That the restrictede weak-type $(1,1)$ property of the ${\cal M}$ implies the almost everywhere
boundary limit result is standard.

The proof of the reverse direction follows the lines in the last section, using of course the asymptotic
estimate $(*)$ for $P$.  In particular, it is straightforward to see that
$$
\lim_{r \ra 0^+} \frac{1}{\sigma(\beta(y,r)} \, \int_{\beta(y,r)} f(t) \, d\sigma(t) = f(y)
$$
if and only if
$$
\lim_{\Gamma_\alpha(y) \ni x \ra y} \int_{\partial \Omega} P(x, t) f(t) \d\sigma(t) = f(x) \, .
$$
Indeed, it is standard to compare the Poisson kernel to a weighted sum of characteristic functions
of balls (the inequality $(*)$ facilitates this comparison), and that gives the result.
That completes the proof.
\endpf
\smallskip \\

\section{Concluding Remarks}

The results of this paper round out the picture of the equivalence of
various key results in basic harmonic analysis.  They will come as
no surprise to the experts, but it useful to have the ideas
recorded in one place, and to know that all these key ideas are logically
equivalent.

\newpage
	     
\null \vspace*{-.3in}

\noindent {\Large \sc References}
\bigskip  \\

\begin{enumerate}

\item[{\bf [APF]}]  L. S. Apfel, Localization properties and
boundary behavior of the Bergman kernel, thesis, 
Washington University in St.\ Louis, 2003.

\item[{\bf [COW]}] R. R. Coifman and G. Weiss, {\it Analyse Harmonique
Non-Commutative Sur Certains Espaces Homogènes; \'{E}tude de Certaines
Int\'{e}grales Singuli\'{e}res}, Springer-Verlag, New York, 1971.

\item[{\bf [GIL]}] J. Gilbert, Niki\u{s}in-Stein
theory and factorization with applications, {\it Harmonic Analysis
in Euclidean Spaces} (Proc.\ Sympos.\ Pure Math., Williams Coll.,
Williamstown, Mass., 1978), Part 2, pp. 233--267, Proc.
Sympos. Pure Math., XXXV, Part, Amer. Math. Soc., Providence,
R.I., 1979.

\item[{\bf [KRA1]}]  S. G. Krantz, {\it Function Theory of Several Complex Variables},
$1^{\rm st}$ ed., John Wiley and Sons, New York, 1982.

\item[{\bf [KRA2]}]  S. G. Krantz, Calculation and estimation of the Poisson kernel,
{\it J. Math.\ Anal.\ Appl.} 302(2005)143--148. 

\item[{\bf [SAW]}]  S. Sawyer, Maximal inequalities of weak type, {\it Ann.\ of Math.} 84(1966), 157--174. 

\item[{\bf [STE1]}]  E. M. Stein, On limits of seqences of operators, {\it Ann.\ of Math.} 74(1961), 140--170. 

\item[{\bf [STE2]}]  E. M. Stein, {\it Singular Integrals and Differentiability Properties of Functions},
Princeton University Press, Princeton, NJ, 1970.

\item[{\bf [STE3]}]  E. M. Stein, {\it Boundary Behavior of Holomorphic Functions
of Several Complex Variables}, Princeton University Press, Princeton, 1972.

\item[{\bf [STW]}]  E. M. Stein and G. Weiss, {\it Introduction to Fourier Analysis on metric Spaces},
Princeton University Press, Princeton, NJ. 1971.

\end{enumerate}

\small

\noindent \begin{quote}
Department of Mathematics \\
Washington University in St. Louis \\
St.\ Louis, Missouri 63130 \\ 
{\tt sk@math.wustl.edu}
\end{quote}

\end{document}